# Hankel determinants of harmonic numbers and related topics

Johann Cigler

Fakultät für Mathematik

Universität Wien

johann.cigler@univie.ac.at

**Abstract.**

We give an overview of known results about Hilbert matrices from the point of view of orthogonal polynomials and compute Hankel determinants of harmonic numbers and related topics.

**Introduction**

We consider harmonic numbers $H_n = \sum_{k=1}^{n} \frac{1}{k}$ and more generally polynomials

$H_n(t,s) = \sum_{k=1}^{n} \frac{st^k}{k+s-1}$ for some $s > 0$. Our aim is the computation of the Hankel determinants

$\det\left(H_{i+j}(t,s)\right)_{i,j=0}^{n}$. In order to make the paper self-contained we also derive some well-known facts about Hilbert matrices $\left(\frac{1}{i+j+1}\right)_{i,j=0}^{n}$ from the point of view of orthogonal polynomials.

I want to thank Tewodros Amdeberhan and Fedor Petrov for providing answers to my questions on MathOverflow, Carsten Schneider for showing me how to use his summation package Sigma, Christian Krattenthaler for alternative proofs of (2.14) and (4.1) and Helmut Prodinger for reference [10].

**1. Some background material**

Let me begin with some well-known background material connected with the approach to Hankel determinants via orthogonal polynomials (cf. e.g. [3],[5],[6],[7]).

Let $(a(n))_{n \geq 0}$ be a sequence of real numbers with $a(0) = 1$.

If no Hankel determinant $\det(a(i+j))_{i,j=0}^{n}$, $n \in \mathbb{N}$, vanishes, then the monic polynomials

$$p(n,x) = \frac{1}{\det(a(i+j))_{i,j=0}^{n-1}} \det \begin{pmatrix} a(0) & a(1) & \cdots & a(n-1) & 1 \\ a(1) & a(2) & \cdots & a(n) & x \\ a(2) & a(3) & \cdots & a(n+1) & x^2 \\ \vdots & & & & \vdots \\ a(n) & a(n+1) & \cdots & a(2n-1) & x^n \end{pmatrix} \quad (1.1)$$

are orthogonal with respect to the linear functional $F$ defined by

$$F(x^n) = a(n). \tag{1.2}$$

This means that $F(p(n,x)p(m,x)) = 0$ if $m \neq n$ and $F(p(n,x)^2) \neq 0$.

Note that $F(p(n,x)x^k) = 0$ for $0 \leq k < n$, because this leads to two identical columns in (1.1).

In particular for $m = 0$ we get

$$F(p(n,x)) = [n = 0]. \tag{1.3}$$

The identities (1.3) also characterize the linear functional $F$.

By Favard's theorem about orthogonal polynomials there exist numbers $s(n), t(n)$ such that

$$p(n,x) = (x - s(n-1))p(n-1,x) - t(n-2)p(n-2,x). \tag{1.4}$$

If on the other hand for given sequences $s(n)$ and $t(n)$ we define numbers $a(n,j)$ by

$$\begin{aligned} a(0,j) &= [j = 0] \\ a(n,0) &= s(0)a(n-1,0) + t(0)a(n-1,1) \\ a(n,j) &= a(n-1,j-1) + s(j)a(n-1,j) + t(j)a(n-1,j+1) \end{aligned} \tag{1.5}$$

then the Hankel determinant $\det(a(i+j,0))_{i,j=0}^{n}$ is given by

$$\det(a(i+j))_{i,j=0}^{n} = \prod_{i=1}^{n} \prod_{j=0}^{i-1} t(j). \tag{1.6}$$

Thus the Hankel determinant $\det(a(i+j,0))_{i,j=0}^{n}$ only depends on the sequence $(t(n))$.

In order to prove this we show first that for all $m, n \in \mathbb{N}$

$$\sum_{k=0}^{n} a(n,k) a(m,k) \prod_{j=0}^{k-1} t(j) = a(m+n, 0). \tag{1.7}$$

We prove this by induction. Since it is true for $n = 0$ and arbitrary $m \in \mathbb{N}$ we assume that it holds for $n-1$ and arbitrary $m \in \mathbb{N}$. Then we get

$$\sum_{k=0}^{n} a(n,k) a(m,k) \prod_{j=0}^{k-1} t(j) = \sum_{k} \left( a(n-1, k-1) + s(k) a(n-1, k) + t(k) a(n-1, k+1) \right) a(m,k) \prod_{j=0}^{k-1} t(j)$$

$$= \sum_{k} a(n-1,k) \left( a(m, k+1) t(k) + s(k) a(m,k) + a(m, k-1) \right) \prod_{j=0}^{k-1} t(j)$$

$$= \sum_{k} a(n-1,k) a(m+1, k) \prod_{j=0}^{k-1} t(j) = a(n+m, 0).$$

Consider now the lower triangular matrix $A_n = (a(i,j))_{i,j=0}^n$ whose diagonal elements are $a(n,n) = 1$ and the diagonal matrix $D_n = \left( [i=j] \prod_{k=0}^{i-1} t(k) \right)_{i,j=0}^n$.

Then (1.7) is equivalent with

$$(a(i+j,0))_{i,j=0}^n = A_n D_n A_n^T \tag{1.8}$$

which implies (1.6).

Thus if we start with $(a(n))$ and guess all $s(n), t(n)$ and $a(n,j)$, then our guesses give correct results if (1.5) holds with $a(n,0) = a(n)$. In this case we also have

$$\sum_{k=0}^n a(n,k) p(k,x) = x^n. \tag{1.9}$$

The polynomials $p(n,x)$ are orthogonal and satisfy

$$F\big((p(n,x))^2\big) = t(0)t(1)\cdots t(n-1). \tag{1.10}$$

This follows by induction because

$$F\big((p(n,x))^2\big) = F\big(x^n p(n,x)\big) = F\big(x^{n-1}(xp(n,x))\big)$$
$$= F\big(x^{n-1}(p(n+1,x) + s(n)p(n,x) + t(n-1)p(n-1,x))\big) = t(n-1) F\big((p(n-1,x))^2\big).$$

Let us also compute the inverse of the Hankel matrices $(a(i+j,0))_{i,j=0}^n$. A simple proof has been given in [3], which I recall in another notation.

Let $T(k) = \dfrac{1}{\prod_{j=0}^{k-1} t(j)}$ and $K_n(x,y) = \sum_{k=0}^n T(k) p(k,x) p(k,y)$ be the so-called kernel polynomial.

Then $K_n(x,y) = \sum_{i,j} b_n(i,j) x^i y^j$ for some numbers $b_n(i,j)$ which satisfy $b_n(i,j) = b_n(j,i)$.

$F\big(x^k K_n(x,y)\big) = \sum_{\ell=0}^n T(k) F\big(x^k p(\ell,x)\big) p(\ell,y)$ is a polynomial of degree $k$ in $y$ because orthogonality implies $F\big(x^k p(i,x)\big) = 0$ for $i > k.$

On the other hand we have

$$F\big(x^k K_n(x,y)\big) = F\left( \sum_{i,j} b_n(i,j) x^{i+k} y^j \right) = \sum_{j=0}^n \sum_{i=0}^n b_n(i,j) a(i+j) y^j.$$

Therefore we get

$$\sum_{i=0}^{n} b_n(i,j)a(i+j) = 0 \text{ for } j > k.$$

For $j = k$ we see that $\sum_{i=0}^{n} b_n(i,k)a(i+k)$ is the coefficient of $y^k$ in

$$\sum_{\ell=0}^{n} T(k)F\left(x^k p(\ell,x)\right)p(\ell,y) \text{ and is therefore equal to}$$

$$T(k)F(x^k p(k,x)) = T(k)F\left((p(k,x))^2\right) = 1.$$

Therefore we get $\left(b_n(i,j)\right)_{i,j=0}^{n} \left(a(i+j)\right)_{i,j=0}^{n} = I_{n+1}$, where $I_{n+1}$ denotes the identity matrix.

Since $b_n(i,j) = b_n(j,i)$ we get

**Proposition 1.1**

Let $\sum_{k=0}^{n} \frac{p(k,x)p(k,y)}{t(0)t(1)\cdots t(k-1)} = \sum_{i,j} b_n(i,j) x^i y^j$. Then the inverse of the Hankel matrix $\left(a(i+j)\right)_{i,j=0}^{n}$ is the matrix

$$\left(\left(a(i+j)\right)_{i,j=0}^{n}\right)^{-1} = \left(b_n(i,j)\right)_{i,j=0}^{n}. \tag{1.11}$$

## 2. Hilbert matrices

Let $a(n) = \dfrac{t^n}{n+1}$ for some $t > 0$. Then $a(n) = F(x^n)$ for the linear functional $F$ on the polynomials defined by $F(p(x)) = \dfrac{1}{t}\int_0^t p(x)dx.$

For $t = 1$ the corresponding Hankel matrix $M_n = \left(\dfrac{1}{i+j+1}\right)_{i,j=0}^{n}$ is the well known Hilbert matrix (cf. e.g. [3],[4]), whose determinant can be directly computed and is given by

$$\det M_n = \frac{1}{\prod_{j=1}^{n}(2j+1)\binom{2j}{j}^2}. \tag{2.1}$$

Let us sketch the above mentioned approach for the special case $F(x^n) = a(n) = \dfrac{t^n}{n+1}$.

Using (1.1) it is easy to guess that $s(n) = \dfrac{t}{2}$ and $t(n) = \dfrac{(n+1)^2 t^2}{4(2n+1)(2n+3)}.$

The first terms of the numbers $a(n,k)$ defined in (1.5) are

$$\begin{pmatrix} 1 & 0 & 0 & 0 & 0 & 0 & 0 & 0 \\ \frac{t}{2} & 1 & 0 & 0 & 0 & 0 & 0 & 0 \\ \frac{t^2}{3} & t & 1 & 0 & 0 & 0 & 0 & 0 \\ \frac{t^3}{4} & \frac{9t^2}{10} & \frac{3t}{2} & 1 & 0 & 0 & 0 & 0 \\ \frac{t^4}{5} & \frac{4t^3}{5} & \frac{12t^2}{7} & 2t & 1 & 0 & 0 & 0 \\ \frac{t^5}{6} & \frac{5t^4}{7} & \frac{25t^3}{14} & \frac{25t^2}{9} & \frac{5t}{2} & 1 & 0 & 0 \\ \frac{t^6}{7} & \frac{9t^5}{14} & \frac{25t^4}{14} & \frac{10t^3}{3} & \frac{45t^2}{11} & 3t & 1 & 0 \\ \frac{t^7}{8} & \frac{7t^6}{12} & \frac{7t^5}{4} & \frac{245t^4}{66} & \frac{245t^3}{44} & \frac{147t^2}{26} & \frac{7t}{2} & 1 \end{pmatrix}$$

This suggests that

$$a(n,k) = \binom{n}{k} \frac{(2k+1)!}{k!} \frac{n!}{(n+k+1)!} t^{n-k}. \tag{2.2}$$

It is now easy to verify that

$$a(n,k) = a(n-1,k-1) + \frac{t}{2} a(n-1,k) + \frac{(k+1)^2 t^2}{4(2k+1)(2k+3)} a(n-1,k+1)$$

and $a(n,0) = \dfrac{t^n}{n+1}$. Therefore all our guesses are correct.

Since

$$\prod_{j=0}^{n-1} t(j) = \frac{t^{2n}}{(2n+1)\binom{2n}{n}^2} \tag{2.3}$$

we see by (1.6) that

$$d(n,t) := \det\left(\frac{t^{i+j}}{i+j+1}\right)_{i,j=0}^{n} = \frac{t^{n^2+n}}{\prod_{j=1}^{n}(2j+1)\binom{2j}{j}^2}. \tag{2.4}$$

Let us illustrate the representation (1.8) for $n = 3$:

$$\begin{pmatrix} 1 & 0 & 0 & 0 \\ \frac{t}{2} & 1 & 0 & 0 \\ \frac{t^2}{3} & t & 1 & 0 \\ \frac{t^3}{4} & \frac{9t^2}{10} & \frac{3t}{2} & 1 \end{pmatrix} \cdot \begin{pmatrix} 1 & 0 & 0 & 0 \\ 0 & \frac{t^2}{12} & 0 & 0 \\ 0 & 0 & \frac{t^4}{180} & 0 \\ 0 & 0 & 0 & \frac{t^6}{2800} \end{pmatrix} \cdot \begin{pmatrix} 1 & \frac{t}{2} & \frac{t^2}{3} & \frac{t^3}{4} \\ 0 & 1 & t & \frac{9t^2}{10} \\ 0 & 0 & 1 & \frac{3t}{2} \\ 0 & 0 & 0 & 1 \end{pmatrix} == \begin{pmatrix} 1 & \frac{t}{2} & \frac{t^2}{3} & \frac{t^3}{4} \\ \frac{t}{2} & \frac{t^2}{3} & \frac{t^3}{4} & \frac{t^4}{5} \\ \frac{t^2}{3} & \frac{t^3}{4} & \frac{t^4}{5} & \frac{t^5}{6} \\ \frac{t^3}{4} & \frac{t^4}{5} & \frac{t^5}{6} & \frac{t^6}{7} \end{pmatrix}$$

The orthogonal polynomials

$$p_n(x,t) = \frac{\prod_{j=1}^{n-1}(2j+1)\binom{2j}{j}^2}{t^{n^2-n}} \det \begin{pmatrix} 1 & \frac{t}{2} & \cdots & \frac{t^{n-1}}{n} & 1 \\ \frac{t}{2} & \frac{t^2}{3} & \cdots & \frac{t^n}{n+1} & x \\ \frac{t^2}{3} & \frac{t^3}{4} & \cdots & \frac{t^{n+1}}{n+2} & x^2 \\ \vdots & & & & \vdots \\ \frac{t^n}{n+1} & \frac{t^{n+1}}{n+2} & \cdots & \frac{t^{2n-1}}{2n} & x^n \end{pmatrix} \qquad (2.5)$$

satisfy

$$p_n(x,t) = \left(x - \frac{t}{2}\right) p_{n-1}(x,t) - \frac{(n-1)^2 t^2}{4(2n-3)(2n-1)} p_{n-2}(x,t). \qquad (2.6)$$

The first terms are

$$\left\{1,\ -\frac{t}{2}+x,\ \frac{t^2}{6}-tx+x^2,\ -\frac{t^3}{20}+\frac{3t^2 x}{5}-\frac{3t x^2}{2}+x^3,\ \frac{t^4}{70}-\frac{2t^3 x}{7}+\frac{9t^2 x^2}{7}-2t x^3+x^4\right\}$$

It is easy to verify that $p_n(0,t) = (-1)^n \dfrac{t^n}{\binom{2n}{n}}$. Therefore we introduce the polynomials

$$P_n(x,t) = \binom{2n}{n} p_n(x,t). \qquad (2.7)$$

The first terms are

$$\left\{1,\ -t+2x,\ t^2-6tx+6x^2,\ -t^3+12t^2 x-30 t x^2+20 x^3,\ t^4-20 t^3 x+90 t^2 x^2-140 t x^3+70 x^4\right\}$$

The coefficients can be found in OEIS [8], A063007, which suggests that

$$P_n(x,t) = \sum_{j=0}^{n} (-t)^{n-j} \binom{n}{j}\binom{n+j}{j} x^j \qquad (2.8)$$

and therefore

$$p_n(x,t) = \sum_{j=0}^{n} (-t)^{n-j} \frac{\binom{n}{j}\binom{n+j}{j}}{\binom{2n}{n}} x^j. \tag{2.9}$$

These formulae can easily be verified by induction.

The polynomials $P_n(x,1)$ are also known as shifted Legendre polynomials.

From (2.6) we conclude that $P_n(x,t)$ satisfies

$$(n+2)P_{n+2}(x,t) - (2x-t)(2n+3)P_{n+1}(x,t) + t^2(n+1)P_n(x,t) = 0 \tag{2.10}$$

with initial values $P_0(x,t) = 1$ and $P_1(x,t) = 2x - t$.

For later use let us mention that (2.10) implies

$$P_n(1,1) = 1. \tag{2.11}$$

Let us also note that

$$\sum_{j=0}^{n} (-t)^{n-j} \binom{n}{j}\binom{n+j}{j} \frac{t^j}{j+1} = F(P_n(x,t)) = \binom{2n}{n} F(p_n(x,t)) = 0 \tag{2.12}$$

for $n > 0$, because orthogonality implies $F(p_n(x,t)) = F(p_n(x,t) p_0(x,t)) = 0$.

Finally we follow Christian Berg [3] to prove the well-known

**Theorem 2.1**

*The inverse matrix of the Hilbert matrix $M_n = \left(\dfrac{1}{i+j+1}\right)_{i,j=0}^{n}$ has integer coefficients.*

*More precisely*

$$M_n^{-1} = \left( (-1)^{i+j}(i+j+1)\binom{n+i+1}{n-j}\binom{n+j+1}{n-i}\binom{i+j}{i}^2 \right)_{i,j=0}^{n}. \tag{2.13}$$

**Proof**

By (1.11) we know that $M_n^{-1} = (b_n(i,j))_{i,j=0}^n$ with

$$\sum_{i,j} b_n(i,j) x^i y^j = \sum_{k=0}^n (2k+1)\binom{2k}{k}^2 p_k(x) p_k(y)$$

$$= \sum_{k=0}^n (2k+1) P_k(x) P_k(y) = \sum_{k=0}^n (2k+1) \sum_{i,j} (-1)^{i+j} \binom{k}{i}\binom{k}{j}\binom{k+i}{i}\binom{k+j}{j} x^i y^j.$$

This implies that $M_n^{-1}$ has integer coefficients.

As in [3] we show that

$$\sum_{k=\max(i,j)}^n \binom{k}{i}\binom{k}{j}\binom{k+i}{i}\binom{k+j}{j}(2k+1) = S(n,i,j) \tag{2.14}$$

with

$$S(n,i,j) = (i+j+1)\binom{n+i+1}{n-j}\binom{n+j+1}{n-i}\binom{i+j}{i}^2$$

by induction in $n$. We can assume that $i \geq j$.

For $n = k = i$ both sides are equal.

We now have

$$S(n+1,i,j) - S(n,i,j) = (i+j+1)\binom{n+i+2}{n+1-j}\binom{n+j+2}{n+1-i}\binom{i+j}{i}^2 - (i+j+1)\binom{n+i+1}{n-j}\binom{n+j+1}{n-i}\binom{i+j}{i}^2$$

$$= (i+j+1)\binom{n+i+1}{n-j}\binom{n+j+1}{n-i}\binom{i+j}{i}^2 \left(\frac{n+i+2}{n+1-j}\frac{n+j+2}{n+1-i} - 1\right)$$

$$= (i+j+1)\binom{n+i+1}{n-j}\binom{n+j+1}{n-i}\binom{i+j}{i}^2 \frac{(i+j+1)(2n+3)}{(n+1-j)(n+1-i)}$$

It remains to show that the last term equals $\binom{n+1}{i}\binom{n+1}{j}\binom{n+1+i}{i}\binom{n+1+j}{j}(2n+3)$,

which is easy to verify.

**Remark 2.2**

Christian Krattenthaler has observed that (2.14) is a special case of the very well-poised hypergeometric summation formula (Lucy Joan Slater [12], Appendix (III.12))

$$_5F_4\left(\begin{array}{c} a,\ \frac{a}{2}+1,\ b,\ c,\ d \\ \frac{a}{2},\ 1+a-b,\ 1+a-c,\ 1+a-d \end{array};1\right) = \frac{\Gamma(1+a-b)\Gamma(1+a-c)\Gamma(1+a-d)\Gamma(1+a-b-c-d)}{\Gamma(1+a)\Gamma(1+a-b-c)\Gamma(1+a-b-d)\Gamma(1+a-c-d)}.$$

Changing $k \to n-k$ in (2.14) gives

$$S(n,i,j) = \frac{(1+2n)(1+n)^{(i)}(1+n)^{(j)}(1+n-i)^{(i)}(1+n-j)^{(j)}}{(i!)^2 (j!)^2} {}_5F_4\left(\begin{array}{c} -1-2n,\ \frac{1}{2}-n,\ i-n,\ j-n,\ 1 \\ -\frac{1}{2}-n,\ -i-n,\ -j-n,\ -i-2n \end{array};1\right),$$

where $x^{(j)} = x(x+1)\cdots(x+j-1)$ denotes a rising factorial. If we set

$\dfrac{\Gamma(-n)}{\Gamma(-m)} = \lim_{\varepsilon \to 0} \dfrac{\Gamma(-n+\varepsilon)}{\Gamma(-m+\varepsilon)} = (-1)^{m-n} \dfrac{m!}{n!}$ for $m,n \in \mathbb{N}$, we get again (2.14).

## 3. Hankel determinants of harmonic numbers

Let $H_n(t) = H_n(t,1) = \sum_{k=1}^{n} \dfrac{t^k}{k}$ and consider the Hankel determinants $D(n,t) := \det\left(H_{i+j}(t)\right)_{i,j=0}^{n}$.

Since $H_0(t) = 0$ the above construction does not work. But by elementary operations we get

$$\det\left(H_{i+j}(t)\right)_{i,j=0}^{n} = (-t)^n \det\begin{pmatrix} 1 & \frac{t}{2} & \cdots & \frac{t^{n-1}}{n} & 0 \\ \frac{t}{2} & \frac{t^2}{3} & \cdots & \frac{t^n}{n+1} & H_1(t) \\ \frac{t^2}{3} & \frac{t^3}{4} & \cdots & \frac{t^{n+1}}{n+2} & H_2(t) \\ \vdots & & & & \vdots \\ \frac{t^n}{n+1} & \frac{t^{n+1}}{n+2} & \cdots & \frac{t^{2n-1}}{2n} & H_n(t) \end{pmatrix}. \quad (3.1)$$

Comparing with (1.1) we see that

$$D(n,t) = (-t)^n d(n-1,t) \sum_{j=0}^{n} (-t)^{n-j} \frac{\binom{n}{j}\binom{n+j}{j}}{\binom{2n}{n}} H_j(t). \quad (3.2)$$

Therefore the computation of the Hankel determinants can be reduced to the computation of the polynomials

$$r(n,t) := \sum_{j=0}^{n} (-t)^{n-j} \binom{n}{j}\binom{n+j}{j} H_j(t). \quad (3.3)$$

**Theorem 3.1**

*Let*

$$r(n,t) = \sum_{j=0}^{n} (-t)^{n-j} \binom{n}{j} \binom{n+j}{j} H_j(t). \tag{3.4}$$

*Then*

$$D(n,t) = \det\left(H_{i+j}(t)\right)_{i,j=0}^{n} = (-1)^n \frac{t^{n^2} r(n,t)}{\binom{2n}{n} \prod_{j=1}^{n-1}(2j+1)\binom{2j}{j}^2}. \tag{3.5}$$

Thus our problem has been reduced to the study of $r(n,t)$. The summation package Sigma by Carsten Schneider [11] gives a computer proof of

**Lemma 3.2**

*The sequence $(r(n,t))$ satisfies the recurrence*

$$nr(n,t) + (t-2)(2n-1)r(n-1,t) + t^2(n-1)r(n-2,t) = 0 \tag{3.6}$$

*with initial values $r(0,t) = 0$ and $r(1,t) = 2t$.*

**Proof**

For the polynomials $f(n,j,t) = (-t)^{n-j}\binom{n}{j}\binom{n+j}{j}H_j(t)$ Sigma provides another set of polynomials $g(n,j,t) = \dfrac{(-1)^{n-j}\binom{n}{j}\binom{n+j}{j}(4n+6)t^{n+2-j}}{n+1-j}\left(\dfrac{(j+j^2+jn)t^j}{(n+1)(n+2)} - \dfrac{j^2 \sum_{1}^{j}\frac{t^\ell}{\ell}}{(n+j-2)}\right)$

such that

$w(n,j,t) := (n+1)t^2 f(n,j,t) + (2n+3)(t-2)f(n+1,j,t) + (n+2)f(n+2,j,t) = g(n,j,t) - g(n,j+1,t)$

for $0 \leq j < n$ and

$w(n,n,t) := (n+1)t^2 f(n,n,t) + (2n+3)(t-2)\bigl(f(n+1,n,t) + f(n+1,n+1,t)\bigr)$
$+ (n+2)\bigl(f(n+2,n,t) + f(n+2,n+1,t) + f(n+2,n+2,t)\bigr) = g(n,n,t) - g(n,0,t).$

By summing over all $0 \leq j \leq n$ we get

$$(n+2)r(n+2,t)+(t-2)(2n+3)r(n+1,t)+t^2(n+1)r(n,t) = \sum_{j=0}^{n} w(n,j)$$

$$= g(n,0,t) - g(n,1,t) + g(n,1,t) - g(n,2,t) + \cdots + g(n,n,t) - g(n,0,t) = 0.$$

Although these results have automatically been found they can be verified by hand.

## 4. The special cases $t = 1$ and $t = 2$.

### 4.1 The Hankel determinants $\det\left(H_{i+j}\right)_{i,j=0}^{n}$

By (3.6) the sequence $r(n,1) = \sum_{j=0}^{n} (-1)^{n-j} \binom{n}{j}\binom{n+j}{j} H_j$ satisfies

$$nr(n,1) - (2n-1)r(n-1,1) + (n-1)r(n-2,1) = 0$$

with initial values $r(0,1) = 0$ and $r(1,1) = 2$.

The sequence $H_n$ satisfies the same recurrence

$$nH_n - (2n-1)H_{n-1}(t) + (n-1)H_{n-2}(t) = 0,$$

because $n(H_n - H_{n-1}) = 1 = (n-1)(H_{n-1} - H_{n-2})$.

The only difference are the initial values $H_0 = 0$ and $H_1 = 1$.

Therefore

$$r(n,1) = 2H_n. \tag{4.1}$$

From (3.2) we deduce

**Theorem 4.1**

$$D(n,1) = \det\left(H_{i+j}\right)_{i,j=0}^{n} = (-1)^n \frac{2H_n}{\binom{2n}{n}\prod_{j=1}^{n-1}(2j+1)\binom{2j}{j}^2}. \tag{4.2}$$

**Remark 4.2**

Perhaps it is instructive to sketch the genesis of this result.

After computing the first terms of the sequence of Hankel determinants

$$\left\{0, -1, \frac{1}{24}, -\frac{11}{129\,600}, \frac{1}{101\,606\,400}, -\frac{137}{2\,016\,379\,008\,000\,000}, \frac{1}{35\,133\,387\,835\,392\,000\,000}, \right.$$

$$\left. -\frac{1}{1\,368\,579\,806\,263\,600\,939\,008\,000\,000}, \frac{761}{658\,299\,967\,151\,148\,396\,655\,182\,662\,860\,800\,000\,000}\right\}$$

I noticed the primes $11, 137, 761, 7129, \cdots$ in the numerators which also appear in the sequence $(H_n)$ whose first terms are

$$0, 1, \frac{3}{2}, \frac{11}{6}, \frac{25}{12}, \frac{137}{60}, \frac{49}{20}, \frac{363}{140}, \frac{761}{280}, \frac{7129}{2520}, \ldots.$$

This observation suggested that $D(n,1)$ is a multiple of $H_n$. The remaining sequence $\frac{D(n,1)}{H_n}$ contains only small primes which indicated that it is a product of factorials and binomials which easily could be found. Thus I had the result but no proof except its reduction to (3.2). Since I could not find anything in the literature I posted my conjecture at MathOverflow, where Fedor Petrov [9] gave the following cute proof of (4.1).

Using $(-1)^k \binom{n+k}{k} = \binom{-n-1}{k}$ gives

$$(-1)^n P_n(x,1) = \sum_{j=0}^{n} \binom{n}{j}\binom{-n-1}{j} x^j = [z^n](1+z)^n (1+xz)^{-n-1}.$$

Since $H_n = \int_0^1 \frac{1-x^n}{1-x} dx$ we get

$$(-1)^n r(n,1) = (-1)^n \int_0^1 \frac{P_n(x,1) - P_n(1,1)}{x-1} dx = [z^n] \int_0^1 \frac{(1+z)^n (1+xz)^{-n-1} - (1+z)^n (1+z)^{-n-1}}{x-1} dx$$

$$= [z^n] \int_0^1 \frac{\left(\frac{1+z}{1+xz}\right)^n \frac{1}{1+xz} - \frac{1}{1+z}}{x-1} dx.$$

Setting $t = \frac{1+z}{1+xz}$ we get $\frac{dt}{dx} = -t \frac{z}{1+xz} = -\frac{t(1-t)}{x-1}$.

Therefore the integral becomes

$$\frac{-1}{1+z} \int_1^{1+z} \frac{1-t^{n+1}}{t(1-t)} dt = \frac{-1}{1+z} \int_1^{1+z} \left(\frac{1}{t} + 1 + t + \cdots + t^{n-1}\right) dt$$

$$= \frac{-\log(1+z) + H_n}{1+z} - \sum_{k=1}^{n} \frac{(1+z)^{k-1}}{k}.$$

Since $-\frac{\log(1-z)}{1-z} = \sum_{k \geq 1} \frac{x^k}{k} \sum_{\ell \geq 0} z^\ell = \sum_{n \geq 0} H_n z^n$ is the generating function of the harmonic numbers we get

$$(-1)^n r(n,1) = [z^n] \frac{-\log(1+z) + H_n}{1+z} = 2(-1)^n H_n$$

and thus (4.1).

After that T. Amdeberhan [1] gave another interesting proof using the Gosper - Zeilberger algorithm and observing that $H_j = [x]\binom{x+j}{j}$.

In a somewhat different terminology he considered the polynomials
$$a_n(x) = \sum_{j=0}^{n}(-1)^j \binom{n}{j}\binom{n+j}{j}\binom{x+j}{j} \text{ which satisfy}$$
$$[x]a_n(x) = \sum_{j=0}^{n}(-1)^j \binom{n}{j}\binom{n+j}{j} H_j = (-1)^n r(n,1).$$

Zeilberger's algorithm gives the recurrence
$$(n+2)^2 a_{n+2}(x) + (2n+3)(2x+1)a_{n+1}(x) - (n+1)^2 a_n(x) = 0.$$

Observing that $[x]x\binom{x+j}{j} = 1$ and $a_n(0) = (-1)^n P_n(0,1) = (-1)^n$ this gives

$$(-1)^{n+2}(n+2)^2 r(n+2,1) + 2(2n+3)(-1)^{n+1} + (-1)^{n+1}(2n+3)r(n+1,1) - (-1)^n (n+1)^2 r(n,1) = 0.$$

By induction we can assume that $r(n,1) = 2H_n$ and $r(n+1,1) = 2H_{n+1}$. Therefore we get

$$(n+2)^2 r(n+2,1) + 2(2n+3)(-1)^{n+1} + (-1)^{n+1}(2n+3)2H_{n+1} + (-1)^{n+1} 2(n+1)^2 H_n = 0.$$

This can be simplified to give $(n+2)^2 r(n+2,1) = 2(n+2)^2 H_{n+2}$ and thus again $r(n,1) = 2H_n$.

The Gosper- Zeilberger algorithm reminded me that Carsten Schneider [11] has a Mathematica package Sigma dealing with multiple sums. I obtained his package and a tutorial by him which finally led to Lemma 3.1.

After posting the first version of this paper Christian Krattenthaler told me another trick to obtain (4.1). It depends on the fact that the derivative
$$\frac{d}{d\varepsilon} \frac{1}{\varepsilon(\varepsilon+1)\cdots(\varepsilon+n-1)} = \frac{-1}{\varepsilon(\varepsilon+1)\cdots(\varepsilon+n-1)}\left(\frac{1}{\varepsilon} + \frac{1}{\varepsilon+1} + \cdots \frac{1}{\varepsilon+n-1}\right)$$

tends to $-\dfrac{H_n}{n!}$ for $\varepsilon \to 1$.

Therefore

$$u(n,\varepsilon) := \sum_{j=0}^{n}(-1)^{n-j} \frac{\binom{n+j}{j} n!}{(n-j)!\varepsilon(\varepsilon+1)\cdots(\varepsilon+n-1)} \text{ satisfies } \left.\frac{du(n,\varepsilon)}{d\varepsilon}\right|_{\varepsilon=1} = -r(n,1).$$

On the other hand by Chu-Vandermonde

$$u(n,\varepsilon) = (-1)^n {}_2F_1\!\left(\begin{array}{c}-n, n+1\\ \varepsilon\end{array}\right) = (-1)^n \frac{(\varepsilon-2)(\varepsilon-3)\cdots(\varepsilon-n-1)}{\varepsilon(\varepsilon+1)\cdots(\varepsilon+n-1)}.$$

This implies $(-1)^{n-1}\dfrac{du(n,\varepsilon)}{d\varepsilon}=\sum_{i=2}^{n+1}\dfrac{u(n,\varepsilon)}{\varepsilon-i}-\sum_{i=0}^{n-1}\dfrac{u(n,\varepsilon)}{\varepsilon+i}$ which converges to $(-1)^{n-1}2H_n$.

Therefore we get again $r(n,1)=2H_n$.

Finally Helmut Prodinger in a comment to the above mentioned Mathoverflow posting provided a reference to his paper [10] which also contains a proof of (4.1).

**4.2 The Hankel determinants** $\det\left(H_{i+j}(2)\right)$

For $t=2$ we get from (3.6)

$$nr(n,2)+4(n-1)r(n-2,2)=0 \tag{4.3}$$

with $r(0,2)=0$ and $r(1,2)=4$.

This gives $r(2n,2)=0$ and $r(2n+1,2)=(-1)^n\dfrac{n!2^{3n+2}}{(2n+1)!!}=(-1)^n\dfrac{(2n+1)!2^{2n+2}}{((2n+1)!!)^2}$.

Therefore we get

**Theorem 4.3**

$$\det\left(H_{i+j}(2)\right)_{i,j=0}^{2n}=0 \tag{4.4}$$

*and*

$$\det\left(H_{i+j}(2)\right)_{i,j=0}^{2n+1}=\dfrac{(-1)^{n+1}2^{4n^2+7n+3}n!}{(2n+1)!(2n+1)!!\prod_{j=1}^{2n+1}\binom{2j}{j}\binom{2j-1}{j}}. \tag{4.5}$$

**Remark 4.4**

T. Amdeberhan [2] has found a direct proof of (4.3).

By (2.10) the polynomials $P_n(x,2)=\sum_{j=0}^{n}(-2)^{n-j}\binom{n}{j}\binom{n+j}{j}x^j$

satisfy

$$(n+2)P_{n+2}(x,2)-2(x-1)(2n+3)P_{n+1}(x,2)+4(n+1)P_n(x,2)=0. \tag{4.6}$$

By (2.12)

$$\int_0^2 P_n(x,2)dx=\int_0^2\sum_{j=0}^{n}(-2)^{n-j}\binom{n}{j}\binom{n+j}{j}x^j=\sum_{j=0}^{n}(-2)^{n-j}\binom{n}{j}\binom{n+j}{j}\dfrac{2^{j+1}}{j+1}=0.$$

On the other hand we have

$$\int_0^2 \frac{P_n(x,2) - P_n(1,2)}{x-1} dx = \sum_{j=0}^n (-2)^{n-j} \binom{n}{j}\binom{n+j}{j} \int_0^2 \frac{x^j - 1}{x-1} dx$$

$$= \sum_{j=0}^n (-2)^{n-j} \binom{n}{j}\binom{n+j}{j} \sum_{\ell=1}^n \frac{2^\ell}{\ell} = r(n,2).$$

Dividing (4.6) by $x-1$ and integrating we get

$(n+2)r(n+2,2) + 4(n+1)r(n,2) = 0$, which is (4.3).

## 5. A more general situation

Following Christian Berg [3] we consider more generally the Hankel matrices

$$\left( \frac{st^{i+j}}{i+j+s} \right)_{i,j=0}^n \quad \text{for } s > 0.$$

Since most proofs can be found in [3] and are a direct generalization of the case $s=1$ we state the results without proof.

Define a linear functional $F$ by

$$F(x^n) = a(n) = \frac{st^n}{n+s}. \tag{5.1}$$

We get

$$s(n) = \frac{2n^2 + (2n-1)s + s^2}{(s+2n-1)(s+2n+1)} t \tag{5.2}$$

and

$$t(n) = \frac{(n+1)^2(n+s)^2 t^2}{(s+2n)(s+2n+1)^2(s+2n+2)}. \tag{5.3}$$

With these values we compute $a(n,k)$ and get

$$a(n,k) = \binom{n}{k} \prod_{j=0}^k \frac{s+k+j}{s+n+j} t^{n-k} \tag{5.4}$$

with $a(n,0) = \frac{st^n}{s+n}$. Therefore all guesses are correct.

This implies that

$$\det\left(\frac{st^{i+j}}{i+j+s}\right)_{i,j=0}^{n} = \frac{s^n t^{n^2+n}}{\prod_{j=1}^{n}(2j+s)\binom{2j+s-1}{j}^2}. \tag{5.5}$$

The monic orthogonal polynomials are given by

$$p_n(x,s) = \frac{1}{\binom{2n+s-1}{n}} \sum_{j=0}^{n}(-t)^{n-j}\binom{n}{j}\binom{n+j+s-1}{n}x^j. \tag{5.6}$$

Here we have

$$F\left((p_n(x,s))^2\right) = t(0)\cdots t(n-1) = \frac{st^{2n}}{(2n+s)\binom{2n+s-1}{n}^2} \tag{5.7}$$

For the inverse matrix we get as in [3]

$$\left(\left(\frac{s}{i+j+s}\right)_{i,j=0}^{n}\right)^{-1} = \left((-1)^{i+j}\frac{i+j+s}{s}\binom{n+i+s}{n-j}\binom{n+j+s}{n-i}\binom{i+j+s-1}{i}\binom{i+j+s-1}{j}\right)_{i,j=0}^{n}. \tag{5.8}$$

The entries of (5.8) can also be written as

$$\frac{(-1)^{i+j}}{s(s+i+j)}\frac{1}{n!^2}\binom{n}{i}\binom{n}{j}(s+i)^{(n+1)}(s+j)^{(n+1)},$$

where $x^{(j)} = x(x+1)\cdots(x+j-1)$ is a rising factorial.

Let now

$$H_n(t,s) = \sum_{k=1}^{n}\frac{st^k}{k+s-1} \tag{5.9}$$

and

$$r(n,t,s) = \sum_{j=0}^{n}(-t)^{n-j}\binom{n}{j}\binom{n+j+s-1}{n}H_j(t,s). \tag{5.10}$$

Then we get in an analogous way as above

$$\det\left(H_{i+j}(t,s)\right)_{i,j=0}^{n} = \frac{(-t)^n}{\binom{2n+s-1}{n}}\det\left(\frac{st^{i+j}}{i+j+s}\right)_{i,j=0}^{n-1}r(n,t,s). \tag{5.11}$$

The summation package Sigma gives the recurrence

$(n+1)(n+s)(2n+3+s)t^2 r(n,t,s)$
$-(2n+2+s)\big((2n+1)(2n+3)+4s(n+1)+s^2 - 2(n+1)^2 t - st(2n+1+s)\big)r(n+1,t,s)$
$+(n+2)(n+1+s)(2n+1+s)r(n+2,t,s) = 0.$

For $t=1$ this reduces to

$(n+1)(n+s)(2n+3+s)r(n,1,s) - (2n+2+s)\big(1+4n+2n^2 + s(2n+3)\big)r(n+1,1,s)$
$+(n+2)(n+1+s)(2n+1+s)r(n+2,1,s) = 0,$

for which again Sigma gives the solution

$$r(n,1,s) = sH_n + H_n(1,s). \tag{5.12}$$

The first terms of $r(n,1,s)$ are

$$\left\{0,\ 1+s,\ \frac{(2+s)(1+3s)}{2(1+s)},\ \frac{(3+s)(4+18s+11s^2)}{6(1+s)(2+s)},\ \frac{(4+s)(18+99s+98s^2+25s^3)}{12(1+s)(2+s)(3+s)}\right\}$$

The first terms of the sequence of Hankel determinants are

$$\left\{0,\ -1,\ \frac{s(1+3s)}{(1+s)^3(2+s)(3+s)},\ -\frac{4s^2(4+18s+11s^2)}{(1+s)^3(2+s)^4(3+s)^2(4+s)^2(5+s)},\right.$$
$$\left.\frac{288s^3(18+99s+98s^2+25s^3)}{(1+s)^3(2+s)^4(3+s)^5(4+s)^3(5+s)^3(6+s)^2(7+s)}\right\}$$

More generally we have

**Theorem 5.1**

$$\det\big(H_{i+j}(1,s)\big)_{i,j=0}^n = \frac{(-1)^n s^{n-1}\big(sH_n + H_n(1,s)\big)}{\binom{2n+s-1}{n}\prod_{j=1}^{n-1}(2j+s)\binom{2j+s-1}{j}^2}. \tag{5.13}$$

**Remark 5.2**

As for Hilbert matrices it would also be interesting to compute the inverses of the Hankel matrices of harmonic numbers. Computer experiments suggest that perhaps the following assertion might be true:

Let $2H_n = \dfrac{U_n}{V_n}$. Then $U_n\big((H_{i+j})_{i,j=0}^n\big)^{-1}$ is an integer-valued matrix.